\documentclass{ws-ijmpc}

\usepackage{ulem}
\usepackage{hyperref}
\usepackage{lipsum,verbatim}
\usepackage{graphicx}
\usepackage{multirow,rotate,subfigure}
\usepackage{amsmath,amssymb}
\usepackage{booktabs}

\usepackage{color}
\definecolor{red}{rgb}{1,0,0}
\definecolor{green}{rgb}{0,1,0}
\definecolor{blue}{rgb}{0,0,1}

\begin{document}

\markboth{Jason A.C.~Gallas}
{Equivalence of orbital equations}

\catchline{}{}{}{}{}

\title{Equivalence among orbital equations of polynomial maps}

\author{Jason A.C.~Gallas} 

\address{
Instituto de Altos Estudos da Para\'\i ba,
  Rua Silvino Lopes 419-2502,\\   
  58039-190 Jo\~ao Pessoa, Brazil,\\
Complexity Sciences Center, 9225 Collins Ave.~1208, 
Surfside FL 33154, USA,\\
  Max-Planck-Institut f\"ur Physik komplexer Systeme,
             01187 Dresden, Germany\\
jason.gallas@gmail.com}

\maketitle

\begin{history}
\received{19 June 2018}
\accepted{4 July 2018}
\centerline{Published 13 August 2018}
\centerline{https://doi.org/10.1142/S0129183118500821}
\end{history}

\begin{abstract}
This paper shows that orbital equations generated by iteration of polynomial
maps do not have necessarily a unique representation.
Remarkably, they may be represented in an infinity of ways, all interconnected
by certain nonlinear transformations.
Five direct and five inverse transformations are established explicitly
between a pair of orbits defined by cyclic quintic polynomials with real roots
and minimum discriminant.
In addition, infinite sequences of transformations generated recursively are
introduced and shown to produce unlimited supplies of equivalent orbital
equations. Such transformations are generic, valid for arbitrary
dynamics governed by algebraic equations of motion.

\end{abstract}

\ccode{PACS Nos.:
      02.70.Wz, 
      02.10.De, 
      03.65.Fd} 


\section{Introduction}
Periodic dynamics of classical oscillators are commonly studied numerically
when the period grows. By contrast, the exact investigation of periodic orbits,
feasible for systems governed by algebraic equations of motion,
is rarely reported because of the inherent theoretical complications.
This paper describes an exploration of exact periodic equations of motion
which arise in the investigation of dynamics generated by polynomial maps.

Many applications, especially in biology and theoretical physics, can be usefully
described by studying the algebraic properties of the equations of motion obtained
by composition of polynomial maps.
Whereas numerical work provides immediate access to dynamical processes,
there is also great merit in exact algebraic work in that it opens the possibility of
uncovering the systematics behind recurring regularities and obtaining valuable
information buried behind the regular self-similar repetition of structures
typically present in phase diagrams.
Of particular interest is the prediction of metric properties of the accumulations of
doubling and adding cascades observed abundantly in applications.
For instance, there are a number of recent surveys concerning the complexities
observed in
accumulations of doubling and adding cascades in laser systems see Ref.~\cite{ato},
in chemistry\cite{field}, in biochemical models\cite{bio1,bio2}, and in the dynamics
of cancer\cite{cancer}.

The rationale for considering equations of motion generated by repeated
iteration is the fact that  equations of higher degree are expected to
comprehend equations of lower degree, so that the solution of higher degree
equations should involve a methodology similar to the one used for equations of
lower degrees. Recall that equations of motion generated by iteration are necessarily
{\sl Abelian equations}, meaning that they can be solved
algebrically\cite{abel,mat}.
Basically, the goal is to explore group properties of Abelian equations generated
by iteration to understand and predict the sprouting of stability in phase diagrams,
to articulate a theoretical framework that could accommodate the insight won
with numerical computations.

This paper reports the discovery of infinite sequences of transformations
interconnecting orbital equations of motion of arbitrary degrees generated by
compositions of i) the quadratic map, ii) the H\'enon map, and
iii) the canonical quartic map, all defined in section \ref{sec:3} below.
Specifically, we show that a certain transformation found to establish the
isomorphism between two totally real cyclic quintic fields of minimum
discriminant is surprisingly not unique and that four additional transformations
exist that establish the same isomorphism.
Moreover, we determine explicitly the five inverse transformations, an important
question not addressed in previous works.
The combined actions on orbital points of the direct, inverse, and composed
transformations is characterized. 
In addition,  generic families of transformations generated recursively are
introduced and use them to produce an unlimited supply of minimum discriminant
isomorphic orbits.
Such transformations are not restricted to the illustrative examples discussed here
but are valid generically for any dynamical systems of
algebraic origin\cite{sil2,kg,ks}. They are important for the study of orbit
proliferation in equivalence classes of equations of motion produced
by discrete maps.

\section{The known equivalence and minimum discriminants}

To motivate the introduction of the generic  transformation chains discussed in
Section \ref{infi}, we start by considering a concrete example.
In 1955, in a pioneering application of computers to algebraic number theory,
Cohn\cite{cohn} produced three tables of irreducible quintic polynomials with integral
coefficients, arranged by increasing order of their discriminant, for the
three {\sl signatures}\cite{ha64b} $(n,\ell)$, namely  $(1,2)$,
$(3,1)$ and $(5,0)$, where $n$ refers to the number of real roots while  $\ell$
refers to the number of pairs of complex roots. 
The cyclic quintics of minimum discriminant $\Delta=14641=11^4$ found by Cohn are
\begin{eqnarray}
  V(x) &=& x^5-\phantom{2}x^4-4x^3+3x^2+3x-1,
             \qquad (\hbox{Vandermonde's quintic})
                                \label{VVV}\\
  G(x) &=& x^5+2x^4-5x^3-2x^2+4x-1.  \quad
                               \label{ggg}
\end{eqnarray}
Cohn manifested his surprise by the fact that his tables contained up to
six distinct polynomials  sharing the  same discriminant and asked whether or not
isodiscriminant polynomials would define the same number field.
Subsequently, Hasse posed the same question, conjecturing the possible isomorphism
between three quintics sharing a factor $47^2$ in their discriminant\cite{ha64,hl69}.
The isomorphism conjectured by Hasse was confirmed by Zassenhaus and
Liang\cite{zl69}, who used a $p$-adic method to find explicit generating
automorphisms of the Hilbert class field over $\mathbb Q(\sqrt{-47})$.
Zassenhaus and Liang were apparently unaware of Cohn's earlier conjecture.

The field isomorphisms conjectured by Cohn were confirmed in 1974 by 
Cartier and Roy\cite{cr74} who, using the same $p$-adic method of Zassenhaus
and Liang\cite{zl69},
reported tables containing explicit polynomial transformations interconnecting
Cohn's isodiscriminant quintics for all three signatures.
In particular, Cartier and Roy found the transformation
\begin{equation}
   T(x) = x^4-4x^2 -x+2,  \label{crt}
\end{equation}
to connect the roots of $V(x)$ to the roots of $G(x)$. The connection is as follows.
Let $x_1, \dots, x_5$  denote the roots of Vandermonde's celebrated quintic $V(x)$.
Then the roots of  $G(x)$ are given by $T(x_1), \dots$, $T(x_5)$.
The inverse transformation, allowing the passage from the roots of $G(x)$
to the roots of $V(x)$ was not considered, neither by Zassenhaus and Liang
nor by Cartier and Roy.
While the intention of the earlier investigatiors was clearly to demonstrate
the conjectured isomorphisms explicitly, there are a number of questions left
open that need to be addressed in the context of the explicit determination
of generating automorphisms of the Hilbert class field.

The specific representation of cyclic extensions of the rational field is
important for a large class of physical problems related to the dynamics and
stability properties of polynomial cycles.
In this context, a complication is that, even for field extensions of relatively
small degrees, to find a suitable representation of cyclic number fields is
far from trivial.
For instance, the literature contains a series of papers\cite{q1,q2,q3}
dedicated to finding suitable forms to represent numbers in a cyclic quartic
extension $K$ of the rational field $\mathbb Q$.
According to Hudson et al.~\cite{q4}, a cyclic quartic extension  $K$
of $\mathbb Q$ may be expressed uniquely in the form
\[ K = {\mathbb Q}\Big(\sqrt{A(D+B\sqrt{D}\,)}\;\Big),  \]
where $A, B, C, D$ are integers such that
i) $A$ is square free and odd,
ii) $D=B^2+C^2$ is square free, $B,C>0$, and
ii) the greatest common divisor of $A$ and $D$ is $1$.
No analogous result is known for quintic or higher order fields.

Before proceeding, we mention that
the determination of number fields having the smallest discriminant
has been a research topic since at least some 120 years.
For, according to Mayer\cite{mayer29}, 
the absolut minimum discriminant for both signatures of cubics was 
determined in 1896 by Furtw\"angler\cite{fu1896}.
Mayer himself obtained the minimum discriminants for the degree $n=4$,
discriminants also considered subsequently by
Godwin\cite{godwin1,godwin2,godwin3} and others\cite{bf89,ford}.
Early references for $n=5$ include the work of Cohn\cite{cohn}
and the thesis by Hunter \cite{hu57}.
They both seem to be among the very first to use computers to investigate
discriminants.
More recent work on $n=5$ was done by Pohst\cite{po75}
and by Takeuchi\cite{ta84}.
Diaz y Diaz\cite{dyd91} reported a table containing
1077 totally real number fields of
degree five having a discriminant less than 2 000 000. 
He finds two nonisomorphic fields of discriminant 1810969, a
prime, and two nonisomorphic fields of discriminant 
$1\,891\,377=3^3\cdot 70051$. All other
number fields in his table are characterized by their
discriminants. Among these fields, three are cyclic and four
have a Galois closure whose Galois group is the dihedral
group $D_5$. The Galois closure for all the other fields 
found  has a Galois group isomorphic to the symmetric group $S_5$,
meaning that the underlying quintics cannot be solved by radicals.
Subsequently, Schwarz, Pohst and Diaz y Diaz\cite{spd94}
reported  the determinantion of all algebraic number fields $F$ of degree $5$
and absolute discriminant less than $2\times10^7$ (totally real fields),
respectively $5\times 10^6$ (other signatures).

\section{The new poly-transformations and their actions}
\label{sec:3}

The polynomial $V(x)$, Eq.~(\ref{VVV}), represents an orbital equation
of motion that is obtained for at least three paradigmatic physical
models, in the so-called {\sl generating partition limit}\cite{gl11}:
the quadratic map $x_{t+1} = 2-x_t^2$,
the H\'enon map $(x,y) \mapsto (2-x^2, y)$, and the canonical quartic
map\cite{jg94,jg95}, namely $x_{t+1} = (x_t^2-2)^2-2$.
For details see, e.g.~Refs.~\cite{res,haa}.
The basic motivation for considering isomorphisms of Vandermonde's
celebrated quintic,
$V(x)$, is to see whether or not it is possible to interconnect distinct
periodic orbits among themselves.
While a general answer to this question does not seem easy, 
interesting extensions and generalizations were obtained that cast
light into a promising theoretical framework to formulate such
interconnections. This is what is reported here.

From a systematic search, we find the following transformations
to provide {\it direct} passages from  $V(x)$ to $G(x)$:
\begin{eqnarray*}
 D_1(x) &=& -{x}^{3}+{x}^{2}+3\,x-2 =
            - (x-2) ({x}^{2}+x-1),\\
 D_2(x) &=& -{x}^{3}+2\,x  = -x ( {x}^{2}-2),\\
 D_3(x) &=& \phantom{-}{x}^{3}-{x}^{2}-2\,x+1,\\
 D_4(x) &=& \phantom{-}{x}^{4}-4\,{x}^{2}-x+2 =
             (x-2) (x+1) ({x}^{2}+x-1) = T(x),\\
 D_5(x) &=& -{x}^{4}+{x}^{3}+4\,{x}^{2}-2\,x-3.
\end{eqnarray*}
Clearly, $D_4(x)$ coincides with $T(x)$, Eq.~(\ref{crt}), previously
found by Cartier and Roy.
In addition, we find the following
{\it inverse} automorphisms for the passage from $G(x)$ to $V(x)$:
\begin{eqnarray*}
 I_1(x) &=& \phantom{-}4\,{x}^{4}+10\,{x}^{3}-15\,{x}^{2}-15\,x+9,\\
 I_2(x) &=& \phantom{-}{x}^{4}+2\,{x}^{3}-5\,{x}^{2}-3\,x+3=
           (x+1) ({x}^{3}+{x}^{2}-6\,x+3),\\
 I_3(x) &=& -2\,{x}^{4}-5\,{x}^{3}+7\,{x}^{2}+7\,x-3 =
            -(x+1)  (2\,{x}^{3}+3\,{x}^{2}-10\,x+3),\\
 I_4(x) &=&-{x}^{4}-2\,{x}^{3}+5\,{x}^{2}+2\,x-3,\\
 I_5(x) &=& -2\,{x}^{4}-5\,{x}^{3}+8\,{x}^{2}+9\,x-5.
\end{eqnarray*}

\begin{table}[!tbh]
  \begin{minipage}[b]{0.48\linewidth}
\tbl{Action of the direct transformations  $D_\ell(x_i)\to y_j$.}
{\begin{tabular}{@{}llllll@{}} \toprule
    & $x_1$ & $x_2$ & $x_3$ & $x_4$ & $x_5$\\
    \midrule
  $D_1(x_j)$: \quad & $y_4$ & $y_1$ & $y_2$ & $y_5$ & $y_3$\\
  $D_2(x_j)$: \quad & $y_5$ & $y_2$ & $y_4$ & $y_3$ & $y_1$\\
  $D_3(x_j)$: \quad & $y_1$ & $y_5$ & $y_3$ & $y_2$ & $y_4$\\
  $D_4(x_j)$: \quad & $y_3$ & $y_4$ & $y_5$ & $y_1$ & $y_2$\\
  $D_5(x_j)$: \quad & $y_2$ & $y_3$ & $y_1$ & $y_4$ & $y_5$\\ \botrule
\end{tabular}\label{ta1}}
\end{minipage}
\begin{minipage}[b]{0.48\linewidth}
\tbl{Action of the inverse transformations $I_\ell(y_i)\to x_j$.}
{\begin{tabular}{@{}llllll@{}} \toprule
    & $y_1$ & $y_2$ & $y_3$ & $y_4$ & $y_5$\\
    \midrule
  $I_1(y_j)$: \quad & $x_2$ & $x_3$ & $x_5$ & $x_1$ & $x_4$\\
  $I_2(y_j)$: \quad & $x_5$ & $x_2$ & $x_4$ & $x_3$ & $x_1$\\
  $I_3(y_j)$: \quad & $x_1$ & $x_4$ & $x_3$ & $x_5$ & $x_2$\\
  $I_4(y_j)$: \quad & $x_4$ & $x_5$ & $x_1$ & $x_2$ & $x_3$\\
  $I_5(y_j)$: \quad & $x_3$ & $x_1$ & $x_2$ & $x_4$ & $x_5$\\ \botrule
\end{tabular}\label{ta2}}
\end{minipage}
\end{table}

\noindent
For every $\ell$, $I_\ell(x)$ is the inverse of $D_\ell(x)$.  Manifestly,
these ten transformations prove non-uniqueness of both sets. They significantly
extend current knowledge concerning generating automorphisms\cite{gy}.
The transformations are not all irreducible, and the degree of the inverses
is always four.

Next, let the roots of these polynomials be so named that
$$\begin{array}{llllll lllll}  
V(x): \ & \ x_1&= -1.68,& \ x_2&=-0.83,& \ x_3&=0.28,& \ x_4&=1.30,& \ x_5&=1.91,\\
G(x): \ & \ y_1&= -3.22,& \ y_2&=-1.08,& \ y_3&=0.37,& \ y_4&=0.54,& \ y_5&=1.39.
\end{array}
$$
Applying $D_i(x)$ and $I_i(x)$ to the roots above produces their rearrangement 
as summarized in Tables \ref{ta1} and \ref{ta2}.
From Table \ref{ta1} one sees that $y_5=D_4(x_3)$, while from Table \ref{ta2}
one gets $x_3=I_4(y_5)$, etc. In other words, the tables show unambiguously
that $I_\ell(x)$ produces arrangements inverse to those produced by $D_\ell(x)$,
for $\ell=1,2,3,4,5$.
In other words, the roots of $G(x)$ are invariant under the compositions
$D_\ell(I_\ell(x))$, $\ell=1,2,3,4,5$ while the roots of $V(x)$ are invariant under
$I_\ell(D_\ell(x))$, $\ell=1,2,3,4,5$.
These compositions are of degrees $12$ or $16$, some of them can be
factored over the rational integers.
Their discriminants are usually given by quite large numbers and may contain
large primes, for instance
{\small
\begin{eqnarray*}
  D_5(I_5(x)) &=& 16x^{16}+160x^{15}+344x^{14}-1208x^{13}-4631x^{12}
  +2996x^{11}\cr
  &&\quad +21362x^{10}-2467x^9-49098x^8
  +1069x^7+60211x^6\cr
  &&\quad -6121x^5-37383x^4+10540x^3+8777x^2-4797x+643,
\end{eqnarray*}}\noindent
whose discriminant is \ \
$ 2^{48} \cdot 19^4  \cdot 97 \cdot 103^4 \cdot 50145997294335406244461$. \ \
In contrast, one also finds pairs of almost identical factors:
{\small
\begin{eqnarray*}
  D_4(I_4(x)) &=& (x^4+2x^3-5x^2-2x+2) \times  (x^4+2x^3-5x^2-2x+5)\; \times\cr
  &&   (x^8+4x^7-6x^6-24x^5+22x^4+30x^3-21x^2-10x+5)
\end{eqnarray*}}\noindent
whose discriminant is \ \
$-2^{20} \cdot 3^8 \cdot 5^{12} \cdot 61 \cdot 71 \cdot 1229 \cdot 7691$.
Additional compositions reveal transformations of ever growing degrees 
that produce root arrangements identical to the ones described and that
display no obvious interconnections among them.
Are $V(x)$ and $G(x)$ the only quintics sharing the minimum discriminant?
This is what we investigate next.

The transformations $D_i(x)$ and $I_j(x)$ may be used to obtain infinite
parameterized representations of irreducible isodiscriminant quintics.
Explicitly, $D_i(x)$ produces  the following set of $n$-parameterized
transformations sharing the discriminant $11^4$, independently
of integer $n$:
$$
{\small
\begin{array}{ll}
  P_n^{(1)}(x) &= {x}^{5}- (5\,n+8 ) {x}^{4}
  + ( 10\,{n}^{2}+32\,n+19) {x}^{3}
  -( 10\,{n}^{3}+48\,{n}^{2}+57\,n+4) {x}^{2}\cr
  &\quad+( 5\,{n}^{4}+32\,{n}^{3}+57\,{n}^{2}+8\,n-32) x
  -{n}^{5}-8\,{n}^{4}-19\,{n}^{3}-4\,{n}^{2}+32\,n+23,\\
  P_n^{(2)}(x) &= {x}^{5}-(5\,n-2) {x}^{4}
  + ( 10\,{n}^{2}-8\,n-5) {x}^{3}
  -( 10\,{n}^{3}-12\,{n}^{2}-15\,n+2) {x}^{2}\cr
  &\quad+( 5\,{n}^{4}-8\,{n}^{3}-15\,{n}^{2}+4\,n+4) x
  -{n}^{5}+2\,{n}^{4}+5\,{n}^{3}-2\,{n}^{2}-4\,n-1,\\
  P_n^{(3)}(x) &= {x}^{5}- \left(5\,n-7 \right) {x}^{4}
    + \left( 10\,{n}^{2}-28\,n+13 \right) {x}^{3}
    - \left( 10\,{n}^{3}-42\,{n}^{2}+39\,n-5 \right) {x}^{2}\cr
&\quad + \left( 5\,{n}^{4}-28\,{n}^{3}+39\,{n}^{2}-10\,n-2
\right) x-{n}^{5}+7\,{n}^{4}-13\,{n}^{3}+5\,{n}^{2}+2\,n-1,\\
P_n^{(4)}(x) &= {x}^{5}- \left(5\,n-12 \right) {x}^{4}
   + \left( 10\,{n}^{2}-48\,n+51 \right) {x}^{3}
   - \left( 10\,{n}^{3}-72\,{n}^{2}+153\,n-96 \right) {x}^{2}\cr
&\quad + \left( 5\,{n}^{4}-48\,{n}^{3}+153\,{n}^{2}-192\,
n+80 \right) x-{n}^{5}+12\,{n}^{4}-51\,{n}^{3}+96\,{n}^{2}-80\,n+23,\\
P_n^{(5)}(x) &={x}^{5}- \left(5\,n+13 \right) {x}^{4}
    + \left( 10\,{n}^{2}+52\,n+61 \right) {x}^{3}
    + \left( 10\,{n}^{3}+78\,{n}^{2}+183\,n+119 \right) {x}^{2}\cr
&\quad+ \left( 5\,{n}^{4}+52\,{n}^{3}+183\,{n}^{2}+238
\,n+70 \right) x-{n}^{5}-13\,{n}^{4}-61\,{n}^{3}-119\,{n}^{2}-70\,n+23.
\end{array}
}
$$
These polynomials are not independent due to linear shifts induced by  $D_i(x)$:
$$
  P_n^{(1)} = P_{n+2}^{(2)},\qquad
  P_n^{(2)} = P_{n+1}^{(3)},\qquad  
  P_n^{(3)} = P_{n+1}^{(4)},\qquad
  P_n^{(4)} = P_{n-5}^{(5)},\qquad
  P_n^{(5)} = P_{n+1}^{(1)}.
$$
They produce the same sequence of isodiscriminant cyclic quintics. In particular,
$$
  G(x) =  P_{-2}^{(1)}(x)
       =  P_0^{(2)}(x) 
       =  P_1^{(3)}(x) 
       =  P_2^{(4)}(x) 
       =  P_{-3}^{(5)}(x).
$$  
Similarly, the transformations $I_j(x)$ lead to irreducible
quintics with discriminant $11^4$:
$${\small
\begin{array}{ll}
P_n^{6}(x) &=  {x}^{5}- \left(5\,n-44 \right) {x}^{4}
    + \left( 10\,{n}^{2}-176\,n+770 \right) {x}^{3}
    - \left( 10\,{n}^{3}-264\,{n}^{2}+2310\,n-6699 \right) {x}^{2}\cr
  &\qquad\ + \left( 5\,{n}^{4}-176\,{n}^{3}+2310\,{n}^{2}-13398\,n
+28974 \right) x\cr
&\qquad\ -{n}^{5}+44\,{n}^{4}-770\,{n}^{3}+6699\,{n}^{2}-28974\,n+49841,\\
  P_n^{(7)}(x) &={x}^{5}- \left(5\,n-14 \right) {x}^{4}
  + \left( 10\,{n}^{2}-56\,n+74 \right) {x}^{3}
  - \left( 10\,{n}^{3}-84\,{n}^{2}+222\,n-183 \right) {x}^{2}\cr
&\qquad\ + \left( 5\,{n}^{4}-56\,{n}^{3}+222\,{n}^{2}-366
  \,n+210 \right) x\cr
  &\qquad\ -{n}^{5}+14\,{n}^{4}-74\,{n}^{3}+183\,{n}^{2}-210\,n+89,\\
  P_n^{(8)}(x) &={x}^{5}- \left(5\,n+16 \right) {x}^{4}
  + \left( 10\,{n}^{2}+64\,n+98 \right) {x}^{3}
  - \left( 10\,{n}^{3}+96\,{n}^{2}+294\,n+285 \right) {x}^{2}\cr
&\qquad\ + \left( 5\,{n}^{4}+64\,{n}^{3}+294\,{n}^{2}+570
  \,n+390 \right) x\cr
  &\qquad\ -{n}^{5}-16\,{n}^{4}-98\,{n}^{3}-285\,{n}^{2}-390\,n-199,\\
  P_n^{(9)}(x) &={x}^{5}- \left(5\,n+16 \right) {x}^{4}
  + \left( 10\,{n}^{2}+64\,n+98 \right) {x}^{3}
  - \left( 10\,{n}^{3}+96\,{n}^{2}+294\,n+285 \right) {x}^{2}\cr
  &\qquad\ + \left( 5\,{n}^{4}+64\,{n}^{3}+294\,{n}^{2}+570,
  \,n+390 \right) x\cr
  &\qquad\ -{n}^{5}-16\,{n}^{4}-98\,{n}^{3}-285\,{n}^{2}-390\,n-199,\\
  P_n^{(10)}(x) &={x}^{5}- \left(5\,n+26 \right) {x}^{4}
  + \left( 10\,{n}^{2}+104\,n +266 \right) {x}^{3}
  - \left( 10\,{n}^{3}+156\,{n}^{2}+798\,n+1337 \right) {x}^{2}\cr
      &\qquad\ + \left( 5\,{n}^{4}+104\,{n}^{3}+798\,{n}^{2}
+2674\,n+3298 \right) x\cr
&\qquad\ -{n}^{5}-26\,{n}^{4}-266\,{n}^{3}-1337\,{n}^{2}-3298\,n-3191.
\end{array}
}$$
Clearly, despite the fact that $I_3(x) \neq I_4(x)$, both transformations
produce identical parameterized forms:  $P_n^{(9)}(x) = P_n^{(8)}(x)$.
Analogously as before:
$$
  P_n^{(6)}  = P_{n-6}^{(7)},\qquad
  P_n^{(7)}  = P_{n-6}^{(9)},\qquad
  P_n^{(8)}  = P_{n+12}^{(6)},\qquad    
  P_n^{(9)}  = P_{n-2}^{(10)},\qquad
  P_n^{(10)} = P_{n+2}^{(8)},
$$
\begin{equation*}
  V(x) = P_{-3}^{(9)}(x)
       = P_{-5}^{(10)}(x)
       = P_{-3}^{(8)}(x)
       = P_{9}^{(6)}(x)
       = P_{3}^{(7)}(x).
\end{equation*}  
The shifts above imply cyclic properties of the determinants defining
discriminants.

\section{Infinite chains of discriminant-preserving transformations}
\label{infi}

There is a richer way of generating infinite families of interrelated
but not trivially connected polynomials sharing the same discriminant,
minimal or not.
Observing that $D_4(x)=T(x)=(x^2-2)^2-2-x$ is a composition of a quadratic
function, we are led to introduce families of recursive algorithms
to generate unbounded sequences of discriminant-preserving polynomials,
quintic or not.

Let $t_0(u)= u^2-\alpha_0$ be a polynomial on an arbitrary variable $u$,
and containing an arbitrary parameter $\alpha_0$.
With  $t_0(u)$, build auxiliary polynomials $\{t_i(u)\}$
\begin{equation}
  t_i(u) = t_{i-1}^2 - \alpha_i, \quad\hbox{ for } \quad  i=1, 2, \dots.
  \label{tt}
\end{equation}
where the several $\alpha_i$ are also chosen arbitrarily.
Manifestly, the number of individual sets $\{t_i(u)\}$ is infinite.
As an {\sl ad hoc} example, we consider the specific sequence obtained by
fixing $\alpha_i = 2$ for all $i$.
In this case, the first few auxiliary polynomials are:
\[
\begin{array}{lll}
  t_0(u) &=& u^2 - 2,\\
  t_1(u) &=&  u^4 -4u^2 + 2,\\
  t_2(u) &=&  u^8 -8u^6 + 20u^4 -16u^2 + 2,\\
  t_3(u) &=&  u^{16} -16u^{14} + 104u^{12} - 352u^{10}
             +660u^8 - 672u^6 +336u^4 - 64u^2 +2.
\end{array}
\]
These polynomials are then used to define a chain of transformations,
namely
\begin{equation}
 \fbox{ $T_i(u)  =  t_i(u) -u$} \quad\hbox{ for } \quad  i=0,1, 2, \dots.
       \label{TTT}
\end{equation}
By construction, $T_1(u)$ coincides with $T(x)$ in Eq.~(\ref{crt}), the
transformation found by  Cartier and Roy.
The transformations $T_i(u)$ preserve discriminants, minimal or not,
for quintics of any signature, cyclic or not.
Generalized  chains may be obtained analogously by iterating more
complicated functions and by allowing $\alpha_i$ to vary as the iteration
proceeds.

When applied to $V(x)$ and $G(x)$, the transformations $T_i(u)$ produce
parameterized families  which split dichotomically into either periodic or
non-periodic sequences of irreducible equivalent quintics.
For instance, representing by $T_i(V)$ the operation of
applying the transformation $T_i(x)$ to the roots of $V(x)$, we obtain
a period-five polynomial cycle interconnecting three polynomials which repeat
mod $5$ indefinitely in the following order:
\[{\small\begin{array}{llllll}
  A(x)=T_0(V),\ &G(x)=T_1(V),\
  &G(x)=T_2(V),\ &A(x)=T_3(V),\ &B(x)=T_4(V),\\
  A(x)=T_5(V),\ &G(x)=T_6(V),\
  &G(x)=T_7(V),\ &A(x)=T_8(V), &B(x)=T_9(V),\ \ \hbox{etc},
\end{array} } \]
where $G(x)$ is defined in Eq.~(\ref{ggg}) and
$$\begin{array}{llllll}
   A(x) &=& x^5 + 2x^4-5x^3 -13x^2 -7x -1, 
              \qquad &\Delta_A=& 11^4,\label{cicA}\\
   B(x) &=& x^5 + 2x^4-16x^3 -24x^2 +48x +32, 
              \qquad &\Delta_B=& 11^4\cdot 2^{20}.\label{cicC}
\end{array}
$$
Clearly, $V(x)$ is not part of the above 5-cycle but leads to it.
It is a sort of {\sl preperiodic} equation of motion, mimicking the known
behavior of preperiodic orbital points\cite{res}.
The above cycling of polynomials implies the existence of an infinity
of additional direct transformations of ever increasing degrees allowing
the passage from $V(x)$ to $G(x)$.
In sharp contrast, the analogous sequence obtained from $T_i(G)$, for
$i=0, 1, 2, \dots$,
produces a non-repeating sequence of quintics.

Many other never-repeating sequences of isomorphic irreducible quintics
sharing similar discriminants may be extracted from additional parameterized
families of totally real cyclic quintics. Four examples are
$${\small
\begin{array}{lllll}
  {\mathbb  A}_n(x) &= {x}^{5} &+ \left( 5\,n-9 \right) {x}^{4}
  + \left( 10\,{n}^{2}-36\,n+28 \right) {x}^{3}
  + \left( 10\,{n}^{3}-54\,{n}^{2}+84\,n-35 \right) {x}^{2}\cr
  & &+ \left( 5\,{n}^{4}-36\,{n}^{3}+84\,{n}^{2}-70\,n+15  \right) x
     +{n}^{5}-9\,{n}^{4}+28\,{n}^{3}-35\,{n}^{2}+15\,n-1,\\
  {\mathbb  B}_n(x) &= {x}^{5} &- \left(5\,n+3 \right) {x}^{4}
    +\left( 10\,{n}^{2}+12\,n-3 \right) {x}^{3}
    -\left( 10\,{n}^{3}+18\,{n}^{2}-9\,n-4 \right) {x}^{2}\cr
  &  &+\left( 5\,{n}^{4}+12\,{n}^{3}-9\,{n}^{2}-8\,n+1 \right)x
  -{n}^{5}-3\,{n}^{4}+3\,{n}^{3}+4\,{n}^{2}-n-1,\\
  {\mathbb  C}_n(x) &={x}^{5} &- \left( 5\,n+15 \right) {x}^{4}
  +( 10\,{n}^{2}+60\,n+35) {x}^{3}
  -( 10\,{n}^{3}+90\,{n}^{2}+105\,n+28) {x}^{2}\cr
  & &+ \left( 5\,{n}^{4}+60\,{n}^{3}+105\,{n}^{2}+56\,n +9 \right) x
     -{n}^{5}-15\,{n}^{4}-35\,{n}^{3}-28\,{n}^{2}-9\,n-1,\\
  {\mathbb  D}_n(x) &= {x}^{5} &- \left(5\,n+18 \right) {x}^{4}
  + \left( 10\,{n}^{2}+72\,n+35 \right) {x}^{3}
  - \left( 10\,{n}^{3}+108\,{n}^{2}+105\,n+16 \right) {x}^{2}\cr
  & &+ \left( 5\,{n}^{4}+72\,{n}^{3}+105\,{n}^{2}+32\,n-2 \right) x
     -{n}^{5}-18\,{n}^{4}-35\,{n}^{3}-16\,{n}^{2}+2\,n+1.
\end{array}
}$$
The discriminant of these quintics does not depend of $n$,
being  $11^4$ for 
${\mathbb  A}_n(x)$, ${\mathbb  B}_n(x)$, and ${\mathbb  C}_n(x)$, and
$11^8$  for ${\mathbb  D}_n(x)$.
Families of isodiscriminant quintics parameterized by more than one parameter
can be also generated with $T_i(u)$ but this will not be pursued here.

To conclude this section,
we observe that the above parametric forms share properties similar to
a celebrated parametrized family of quintics, found by Emma Lehmer\cite{e88,sn09}
to provide connections between the so-called Gaussian period equations
and cyclic units. For instance, in the notation of Butler and McKay\cite{bm83},
their common Galois group is 5T1.
It is tempting to conjecture that the parametric forms reported here might
also be linear combinations of Gaussian periods, providing interconnections
between Gaussian period equations and cyclic units.
This, however, remains to be ascertained.

\section{Conclusions and outlook}

This paper reported direct and inverse transformations showing
that in addition to a transformation found by Cartier and Roy,
Cohn's conjectured isomorphism among minimum discriminant cyclic
quintics can be established by nine new transformations,
Therefore, altogether there are five direct and five inverse transformations
which, when combined, reveal how orbital points are rearranged
cyclically under the transformations.

An infinite chain of transformations was introduced and used to show that,
apart from Cohn's pair of quintics, there is an apparently unbounded
quantity of isomorphic quintics sharing the same field discriminant and
group properties.
The chain of transformations is generated by a simple recurrence relation,
Eq.~(\ref{TTT}),
and is valid to transform arbitrary equations of motion, of any degree,
for systems governed by discrete maps. The chain introduced here may be
naturally extended by replacing the quadratic functions underlying
Eq.~(\ref{TTT}) by any arbitrary set of functions.
In other words, the chain of transformations provides an effective tool to
study transformation properties of nonlinear equations of mathematical physics
so popular nowadays in practical applications.
We hope to return to this in the future.

To conclude, we remark that the aforementioned question raised by Hasse
was answered only in part by Zassenhaus and Liang.
It remains to be determined whether or not there are additional analogous
transformations interconnecting Hasse's triplet of quintics, which are
{\sl not cyclic}  and have complex roots.
Explicit expressions for generating automorphisms play a significant
role in the study of Galois-theoretic aspects of iterated maps
and Abelian groups\cite{a1,a2,a3,a4}.

\medskip\medskip

\section*{Acknowledgments}
This work was supported by the
Max-Planck Institute for the Physics of Complex Systems, Dresden,
in the framework of the Advanced Study Group on Optical Rare Events.
The author was also supported by CNPq, Brazil.






\begin{thebibliography}{00}
\bibitem{ato} J.A.C.~Gallas,
  {Spiking systematics in some CO$_2$ laser models},
  {\it Adv. Atom. Molec. Opt. Phys.} {\bf 65},  127--191  (2016).

\bibitem{field} R.J.~Field, 
  Chaos in the Belousov Zhabotinsky reaction,
  {\it Mod. Phys. Lett. B} {\bf 29}, 1530015 (2015).

\bibitem{bio1} J.G.~Freire, M.R.~Gallas, and J.A.C.~Gallas,
  Chaos-free oscillations,
  {\it Europhys. Lett.} {\bf 118}, 38003 (2017).
  
\bibitem{bio2} M.R.~Gallas and J.A.C.~Gallas,
  Nested arithmetic progressions of oscillatory phases in Olsen's enzyme
  reaction model,
  {\it Chaos} {\bf 25}, 064603 (2015).

\bibitem{cancer} M.R.~Gallas, M.R.~Gallas, and J.A.C.~Gallas,
  Distribution of chaos and periodic spikes in a three-cell population model of cancer,
  {\it  European Phys. J. Special Topics} {\bf 223}, 2131--2144 (2014).

  
\bibitem{abel} N.H.~Abel,
  M\'emoire sur une classe particuli\`ere d'\'equations r\'esolubles
  alg\'ebriquement,
  {\it J. reine Angew. Math.} {\bf 4}, 131--156 (1829).
  For a Portuguese translation of this paper see
  {\it Bol. Soc. Portuguesa de Mat.} {\bf 47} 1-15, 17-48 (2002).

\bibitem{mat} G.B.~Mathews, Algebraic Equations
  (Cambridge University Press, Cambridge, 1907).
  
\bibitem{sil2}  J.H. Silverman,
  The Arithmetic of Dynamical Systems
   (Springer, New York, 2007).

\bibitem{kg} K. Gatermann,
      Computer Algebra Methods for Equivariant Dynamical Systems
      (Springer, Berlin, 2000).

\bibitem{ks} K. Schmidt,
   Dynamical Systems of Algebraic Origin
   (Birkh\"auser, Basel, 1995).

\bibitem{cohn} {H.~Cohn}, 
   { A numerical study of quintics of small discriminant},
   {\it Commun. Pure Appl. Math.} {\bf 8},  377--386 (1955).

\bibitem{ha64b} {H.~Hasse},
  Vorlesungen \"uber Zahlentheorie, Zweite Auflage
  (Springer, Berlin, 1964).
   
\bibitem{ha64} {H.~Hasse}, 
  { \"Uber den Klassenk\"orper zum quadratischen Zahlk\"orper mit der
  Diskriminante $-47$},
  {\it Acta Arith.} {\bf 9},  419--434 (1964),

\bibitem{hl69} {H.~Hasse}  and {J.~Liang},   
   { \"Uber den Klassenk\"orper zum quadratischen Zahlk\"orper mit der
  Diskriminante $-47$ (Fortsetzung)},
  {\it Acta Arith.} {\bf 16},  89--97 (1969).

\bibitem{zl69} {H.~Zassenhaus} and {J.~Liang},   
  { On a problem of Hasse},
  {\it Math. Comp.} {\bf 23},  515--519 (1969).
  
\bibitem{cr74} {P.~Cartier} and {Y.~Roy},
   { On the enumeration of quintic fields with small discriminants},
   {\it J. reine angew. Math.} {\bf 268/269},   213--215  (1974).

\bibitem{q1} K.~Hardy, R.H.~Hudson, D.~Richman, K.S.~Williams, and N.M.~Holtz,
  {\it Calculation of class numbers of imaginary cyclic quartic fields},
  Carleton-Ottawa Math. Lecture Note Series {\bf 7}, July 1986, 201 pp.

\bibitem{q2} H.~Edgar and B.~Peterson, Some contributions to the theory of
  cyclic quartic extensions of the rationals,
  {\it J. Number Theory} {\bf 12},  77--83  (1980).

\bibitem{q3} Z.~Xianke,
  Cyclic quartic fields and genus theory of their subfields,
  {\it J. Number Theory} {\bf 18},  350--355  (1984).

\bibitem{q4} R.H.~Hudson and K.S.~Williams,
  The integers of a cyclic quartic field,
  {\it Rocky Mountain J. Math.} {\bf 20},  145--150  (1990).



\bibitem{mayer29} J.~Mayer,
   Die absolut kleinsten Diskriminant der biquadratischen
   Zahlk\"orper,
   {\it Sitzungsber. Akad. Wissens. Wien, math.-naturw. Kl., Abt.~IIa,}
   {\bf 138}, Heft 9/10, 733--742 (1929).

\bibitem{fu1896} Ph.~Furtw\"angler,
   Zur Theorie der in Linearfaktoren zerlegbaren, ganzzahligen
   tern\"aren kubischen Formen,  G\"ottingen, 1896. 
   It was not possible for me to locate this reference.

\bibitem{godwin1} H.J.~Godwin,
   { Real quartic fields with small discriminants},
   {\it J. London Math. Soc.} {\bf 31}, 478--485 (1956).

\bibitem{godwin2} H.J.~Godwin,
   {On quartic fields of signature one with small
        discriminants},
   {\it Quart. J. Math. (Oxford)} {\bf 8}, 214--222 (1956).

\bibitem{godwin3} H.J.~Godwin,
   {On totally complex quartic fields with small discriminants},
   {\it Proc. Cambridge Phil. Soc.} {\bf 53},  1--4 (1957).

\bibitem{bf89} J.~Buchmann and D.~Ford,
   {On the computation of totally real quartic fields of small
        discriminant},
   {\it Math. Comp.} {\bf 52},  161--174 (1989).

\bibitem{ford} D.~Ford,
  {Enumeration of totally complex quartic fields of small
  discriminants}, in ``Computational Number Theory'',
  edited by A.~Peth\"o, M.E. Pohst, H.C. Williams and H.G. Zimmer,
  Walter de Gruyter, Berlin, 1991.

\bibitem{hu57} J.~Hunter,
     The minimum discriminant of quintic fields,
     {\it Proc. Glasgow Math. Assoc.} {\bf 3}, 57--67 (1957).

\bibitem{po75} M.~Pohst, 
   {Berechnung kleiner Diskriminanten total reeller algebraischer
        zahlk\"orper},
   {\it J. reine angew. Math.} {\bf 278/279}, 278--300 (1974).

\bibitem{ta84} K.~Takeuchi,
   Totally real algebraic number fields of degree 5 and 6
   with small discriminants,
   {\it Saitama Math. J.} {\bf 2}, 21--32 (1984).
   
\bibitem{dyd91} F.~Diaz y Diaz,
   A table of totally real quintic number fields,
   {\it Math. Comp.} {\bf 56} (1991), 801--808.
   table supplement: {\bf 56}, S1--S12 (1991).

\bibitem{spd94} A.~Schwarz, M.~Pohst and F.~Diaz y Diaz,
    {A table of quintic number fields}, 
    {\it Math. Comp.} {\bf 63}, 361--376 (1994).

  
\bibitem{gl11} R.~Gilmore and M.~Lefranc,  The Topology of Chaos,
  2nd Edition (Wiley-VCH, Weinheim, 2011).  

  
\bibitem{jg94} J.A.C.~Gallas, 
     Dissecting shrimps: results for some one-dimensional physical systems,
            {\it  Physica A} {\bf 202},   196--223  (1994).

\bibitem{jg95}  J.A.C.~Gallas, 
     Structure of the parameter space of a ring cavity,
     {\it Appl. Phys. B Supplement} {\bf 60},  S203--S213  (1995).

\bibitem{res} J.A.C.~Gallas,
  Method for extracting arbitrarily large orbital equations of
  the Pincherle map,
  {\it Results in Physics} {\bf 6},  561--567 (2016).
  
\bibitem{haa}  J. Argyris, G. Faust, M. Haase, R. Friedrich,
      An Exploration of Dynamical Systems and Chaos, 2nd ed.
      (Springer, New York, 2015).

     
\bibitem{gy} J.H.~Evertse and K.~Gyory,
  Discriminant Equations in Diophantine Number Theory
  (Cambridge University Press, Cambridge, 2017).

\bibitem{e88} E.~Lehmer,
  Connection between Gaussian periods and cyclic units,
  {\it Math. Comp.} {\bf 50}, 535--541 (1988).

\bibitem{sn09} S.~Nakano,
  A family of quintic cyclic fields with even class number parameterized by
  rational points on an elliptic curve,
  {\it J. Number Theory} {\bf 219}, 2943--2951 (2009).

\bibitem{bm83} G.~Butler and J.~McKay, 
  The transitive groups of degree up to eleven,
  {\it Commun. Algebra} {\bf 11}, 863--911  (1983).
  
\bibitem{a1} B.~Eick, The automorphism group of a finitely generated
  virtually abelian group,
  Groups Complexity Cryptology {\bf 8}  (2016) 35--45.

\bibitem{a2} C.~Carstensen, B.~Fine, and G.~Rosenberger,
  Abstract Algebra, Applications to Galois Theory, Algebraic Geometry
  and Cryptography
  (de Gruyter, Berlin, 2011).

\bibitem{a3} J.J.~Cannon and D.F.~Holt, Automorphism group computation
  and isomorphism testing in finite groups,
    {\it J. Symbolic Comput.} {\bf 35},  241--267   (2003).
  
\bibitem{a4} G.~Baumslag, B.~Fine, M.~Kreuzer, and G.~Rosenberger,
 A Course in Mathematical Cryptography    
  (de Gruyter, Berlin, 2015).

\end{thebibliography}
\end{document}